\documentclass[reqno,11pt]{amsart}
\usepackage{amsmath,amsthm}
\usepackage{amsfonts,amssymb,epsfig}
\usepackage[vcentermath,enableskew]{youngtab-op}

\newtheorem{theorem}{Theorem}
\newtheorem{lemma}[theorem]{Lemma}
\newtheorem{proposition}[theorem]{Proposition}
\newtheorem{corollary}[theorem]{Corollary}

\newtheorem{conjecture}{Conjecture}

\def \tr{{\rm tr }}

\newcommand{\bb}{\mathbb}

\def \R {{\bb R}}

\newcommand{\negone}{-1}
\newcommand{\negtwo}{-2}
\newcommand{\T}{{\top}}
\newcommand{\RW}{{{RW}}}

\Yinterspace{2mm}
\Yboxdim{17pt}

\title[Interlacings and the interchange process]%
{Interlacings for random walks on weighted graphs\\
and the interchange process}
\author{A. B. Dieker}
\address{Georgia Institute of Technology, Atlanta GA 30332, USA}
\email{ton.dieker@isye.gatech.edu}

\begin{document}
\begin{abstract}
We study Aldous' conjecture that the spectral gap 
of the interchange process on a weighted undirected graph equals the spectral gap 
of the random walk on this graph.
We present a conjecture in the form of an inequality, 
and prove that this inequality implies Aldous' conjecture
by combining an interlacing result for Laplacians of random walks on weighted graphs
with representation theory.
We prove the conjectured inequality for several important instances.
As an application of the developed theory, we prove Aldous' conjecture for a large class of weighted graphs,
which includes all wheel graphs, all graphs with four vertices, certain nonplanar graphs, certain graphs with several weighted cycles of arbitrary length, as well as all trees.

Caputo, Liggett, and Richthammer have recently resolved Aldous' conjecture,
after independently and simultaneously discovering the key ideas developed in the present paper.
\end{abstract}

\maketitle

\section{Introduction}
This paper studies a fundamental question arising from the theory of card shuffling,
where the evolution of card positions is typically modeled by a Markov chain 
on the space of permutations on the set of cards.
In this paper, we investigate a continuous-time Markov chain in which
the cards at positions $i$ and $j$ are interchanged at rate $\alpha_{ij}$.
Interchange rates may be zero if cards at the corresponding positions cannot be interchanged.
This Markov chain is known as the {\em interchange process}.
Another continuous-time Markov chain arises as the position of an arbitrary but 
fixed card in a deck which evolves according to the interchange process.
This Markov chain is known as the {\em random-walk process}.

A key question is how long it takes for a deck of cards to be well-shuffled in some sense,
and an important quantity in addressing questions of this form is the {\em spectral gap}.
Assuming the interchange rates $\alpha_{ij}$ are chosen so that both the interchange process and the 
random-walk process are irreducible,
the spectral gap is defined as the negative of the second 
largest eigenvalue of the intensity matrix of the interchange process.
A conjecture of Aldous and Diaconis from 1992, often referred to as Aldous' conjecture in the literature,
says that the spectral gap of the interchange process is {\em exactly} the same as the spectral
gap of the random-walk process.
This conjecture,
which is also listed in the open problem section of the recent book by Levin~{\em et al.}~\cite[Sec.~23.3]{MR2466937}, 
is the topic of the present paper.
(Strictly speaking, the original conjecture is more restrictive than in the above discussion, since it only allows $\alpha_{ij}$ to take values in $\{0,1\}$).

It is customary to think of the interchange and random-walk processes in terms of an undirected, connected,
weighted graph $G$.
Each vertex of $G$ has a label, and each edge $(i,j)$ of $G$
has an associated Poisson process with intensity $\alpha_{ij}\ge 0$.
The Poisson processes on different edges are stochastically independent.
At each Poisson epoch corresponding to edge $(i,j)$, the labels at vertices $i$ and $j$ are interchanged.
The interchange process records the positions of all labels in the graph, while 
the random-walk process only records the position of a given label.

\medskip
Aldous' conjecture has attracted the attention of many researchers over the past decades, 
but all existing results rely on some special structure on the weights $\alpha_{ij}$.
These results can roughly be categorized according to their proof methods:
induction on the number of vertices in the graph \cite{MR1419872,MR2415139,starrconomos}
or representation theory~\cite{cesi2,cesi,MR964069,MR626813,MR770635}.
The current paper effectively combines these two approaches, 
and might serve as a first step towards 
proving Aldous' conjecture in full generality.

The main idea behind the combination of mathematical 
induction and representation theory can be summarized as follows.
In the induction step, a new vertex is attached to a graph for which it is known that the
conjecture holds.
If there is only one new edge incident to the new vertex, then
standard eigenvalue bounds can be employed which 
imply that the conjecture holds for the new graph \cite{MR1419872}.
In the general case where several edges are incident to the new vertex, however,
the main technical obstacle has been that the addition of this vertex may significantly impact the spectrum of the 
resulting random walk and interchange process. 
This difficulty can be overcome with representation theory.
In fact, we shall argue that the following conjecture suitably controls the changes to the spectrum
if the new vertex is of degree $k-1$. An empty sum should be interpreted as zero, $\mathcal S_k$
is the symmetric group on $k$ letters, and $(ij)\in \mathcal S_k$ stands for the transposition
of $i$ and $j$.

\begin{conjecture}
\label{conj}
Given any $k\ge 2$, 
the following holds for any function $g:\mathcal S_k\to\R$ and 
any nonnegative $\gamma_{1},\ldots,\gamma_{k-1}$:
\[
\sum_{\sigma\in \mathcal S_k} \sum_{i=1}^{k-1} \gamma_{i} 
[g(\sigma) - g((ik)\sigma)]^2\ge
\sum_{\sigma\in \mathcal S_k}  \sum_{1\le i<j\le k-1} \frac{\gamma_{i}\gamma_{j}}
{\gamma_{1}+\ldots+\gamma_{k-1}} [g(\sigma) - g((ij)\sigma)]^2.
\]
\end{conjecture}
This inequality can be interpreted as a comparison of two Dirichlet forms,
with the left-hand side corresponding to the interchange dynamics on a (weighted) `star' graph with center $k$
and the right-hand side corresponding to the interchange dynamics on a (weighted) special complete graph with 
an isolated vertex $k$.

The contributions of this paper are threefold.
First, we show that Conjecture~\ref{conj} implies Aldous' conjecture. One of the key ingredients is an  interlacing result for Laplacians of random walks on weighted graphs, which appears to be new.
Second, we give a proof of Conjecture~\ref{conj} for $k\le 4$ as well as a proof for general $k$ 
when $\gamma_{1}=\ldots=\gamma_{k-1}$.
Third, as an application of the developed theory, we prove 
that Aldous' conjecture holds for a large family of weighted graphs that only rely on Conjecture~\ref{conj}
for $k\le 4$. 
This class includes all wheel graphs, all weighted graphs with four vertices, certain graphs with weighted cycles of different lengths, certain nonplanar graphs, as well as all trees (for which the conjecture is already known to hold). 
It is the first time results for such general weighted graphs are obtained,
illustrating the power of knowing that Conjecture~\ref{conj} holds even for small values of $k$.

Throughout, matrix inequalities of the form $A\le B$ should
be interpreted as $A-B$ being negative semidefinite.
All vectors in this paper should be interpreted as column vectors, and we use the symbol ${}^{\T}$ for vector or matrix transpose.
The $\ell\times\ell$ identity matrix is denoted by $I_\ell$.
We multiply permutations from right to left, so $\sigma'\sigma$ is the permutation obtained by
first applying $\sigma$ and then $\sigma'$.

This paper is organized as follows. 
Sections~\ref{sec:interlacing}--\ref{sec:conjstronger} focus on proving that Conjecture~\ref{conj} implies Aldous' conjecture. The main technical tools are the aforementioned interlacing for Laplacians of random walks,
which is discussed in Section~\ref{sec:interlacing}, and a representation-theoretic view of Conjecture~\ref{conj},
which is the topic of Section~\ref{sec:representation}.
These tools are tied together in Section~\ref{sec:conjstronger}, which contains the main argument of
the proof that Conjecture~\ref{conj} implies Aldous' conjecture.
Sections~\ref{sec:proofconjS4} and \ref{sec:proofconjJM} prove Conjecture~\ref{conj} in special cases:
Section~\ref{sec:proofconjS4} focuses on $k\le 4$ with general nonnegative $\gamma_i$,
while Section~\ref{sec:proofconjJM} deals with general $k$ but identical $\gamma_i$.
In Section~\ref{sec:application}, we present a class of graphs for which we can prove Aldous' conjecture
using the newly developed methodology.
A discussion concludes this paper, and two appendices give background on
representation theory.

{\bf A postscript; independent work of Caputo, Liggett, and Richthammer.}
Aldous' conjecture was one of three problems targeted by an international team 
of researchers at the Markov Chain Working Group in June 2009, held at Georgia Tech.
I presented this paper at that meeting, and Pietro Caputo presented a joint work 
with Thomas Liggett and Thomas Richthammer.
Both teams had posted their work on arxiv.org 
at the beginning of the meeting \cite{caputov1,dieker:interchangev1}.
Although the papers were written from different perspectives,
we had independently arrived at the same proof outline:
both works propose the updating rule (\ref{eq:defalphaprime}) below
and formulate Conjecture~\ref{conj}.

Only days after the working group meeting, Caputo {\em et al.}~were able to give a full proof of
Conjecture~\ref{conj}. It can be found in \cite{caputo}. 
I currently do not know if it is possible
to give a different proof of Conjecture~\ref{conj} using representation theory, i.e., to complete the 
approach taken here.

The present article is an unmodified version of my `working group' preprint \cite{dieker:interchangev1},
with some typos corrected and some arguments clarified.

\section{Interlacings for Laplacians of random walks on weighted graphs}
\label{sec:interlacing}
In this section, we state and prove an interlacing result for the weighted random-walk process
on a given weighted graph $G$ with $n$ vertices. 
For other interlacing results and illustrations of the technique, 
we refer to Godsil and Royle~\cite[Ch.~9]{MR1829620} or (in a slightly different setting) 
the recent paper by Butler~\cite{MR2320563}.

Let $\alpha_{ij}\ge 0$ be the weight of edge $(i,j)$ in $G$, $i\neq j$.
We simply write $\alpha$ for the collection of edge weights $\{\alpha_{ij}\}$.
We also write $\{e_i\}$ for the standard basis in $\R^n$, and define $w_{ij}=e_i-e_j$.
The Laplacian of $G$ is defined through
\[
\mathcal L_n^\RW(\alpha)(i,j)=
\begin{cases}
-\alpha_{ij} & \text{if $i\neq j$}\\
\sum_{k=1}^{j-1} \alpha_{kj}+
\sum_{k=j+1}^n \alpha_{jk} & \text{if $i=j$}
\end{cases}
\]
for $1\le i,j\le n$, 
and can thus be written in matrix form as
\[
\mathcal L_n^\RW(\alpha) = \sum_{1\le i<j\le n} \alpha_{ij} w_{ij} w_{ij}^\T.
\]
The superscript `RW' is meant to stress that this Laplacian is the negative of the intensity
matrix of the random-walk process defined in the introduction.
Consider edge weights $\alpha_{ij}'$ given by, for $1\le i,j\le n-1$,
\begin{equation}
\label{eq:defalphaprime}
\alpha_{ij}'= \alpha_{ij} + \frac{\alpha_{in}\alpha_{jn}}{\alpha_{1n} + \ldots + \alpha_{n-1,n}},
\end{equation}
while $\alpha_{in}'=0$ for $i<n$.
We abuse notation by writing $\alpha'$ for the restriction to $i,j\le n-1$.
We write $\mu_1^n(\alpha)\le \ldots\le \mu_n^n(\alpha)$ for the eigenvalues of 
$\mathcal L_n^\RW(\alpha)$.
The edge weights $\alpha'_{ij}$ assume a particularly simple form if $\alpha_{in}=0$ for all $i$ except possibly for one: then $\alpha'_{ij}=\alpha_{ij}$ for $i,j\le n-1$.

Note that $\mathcal L_n^\RW(\alpha')$ is the Laplacian of a graph with an isolated vertex $n$.

\begin{proposition}
\label{prop:interlacing}
If $\sum_{i=1}^{n-1}\alpha_{in}>0$, then we have 
\[
\mathcal L_n^\RW(\alpha) = \mathcal L_n^\RW(\alpha') + \frac1{\sum_{i=1}^{n-1}\alpha_{in}}
\left(\sum_{i=1}^{n-1} \alpha_{in} w_{in}\right) \left(\sum_{i=1}^{n-1} \alpha_{in} w_{in}\right)^\T.
\]
In particular, the eigenvalues of $\mathcal L_n^\RW(\alpha)$ and $\mathcal L_n^\RW(\alpha')$ interlace, i.e.,
\[
\mu_1^n(\alpha')\le \mu_1^n(\alpha)\le \mu_2^n(\alpha')\le \ldots \le \mu_n^n(\alpha')\le \mu_n^n(\alpha).
\]
Moreover, for $1\le k\le n$, 
\[
\mu_k^n(\alpha)-\mu_k^n(\alpha')
\le \frac{2\sum_{1\le i\le j \le n-1} \alpha_{in}\alpha_{jn}}{\sum_{i=1}^{n-1} \alpha_{in}}.
\]
\end{proposition}
\proof
First observe that 
\[
\mathcal L_n^\RW(\alpha) - \mathcal L_n^\RW(\alpha') = \sum_{i=1}^{n-1} \alpha_{in} w_{in} w_{in}^\T
-\sum_{1\le i<j\le n-1} \frac{\alpha_{in}\alpha_{jn}}{\alpha_{1n}+\ldots+\alpha_{n-1,n}} w_{ij} w_{ij}^\T.
\]
For the last sum, we note that
\begin{eqnarray*}
\sum_{1\le i<j\le n-1} \alpha_{in}\alpha_{jn} w_{ij} w_{ij}^\T &=&
\sum_{1\le i<j\le n-1} \alpha_{in}\alpha_{jn} (w_{in}-w_{jn}) (w_{in}-w_{jn})^\T\\
&=& (\alpha_{1n}+\ldots+\alpha_{n-1,n}) \sum_{i=1}^{n-1} \alpha_{in} w_{in} w_{in}^\T
- \sum_{1\le i,j\le n-1} \alpha_{in}\alpha_{jn} w_{in} w_{jn}^\T.
\end{eqnarray*}
On combining the preceding two displays, we conclude that
\[
\mathcal L_n^\RW(\alpha)-\mathcal L_n^\RW(\alpha') = 
\sum_{1\le i,j\le n-1} \frac{\alpha_{in}\alpha_{jn}}{\alpha_{1n}+\ldots+\alpha_{n-1,n}} w_{in} w_{jn}^\T,
\]
and we have therefore proven the first claim.

The interlacing property follows from standard results in linear algebra on
the spectrum of rank-1 perturbations of symmetric matrices, e.g., \cite[Thm.~4.3.4]{MR1084815}.
After noting that 
\[
\frac1{\sum_{i=1}^{n-1}\alpha_{in}}
\left(\sum_{i=1}^{n-1} \alpha_{in} w_{in}\right)^\T \left(\sum_{i=1}^{n-1} \alpha_{in} w_{in}\right) = \frac{2\sum_{1\le i\le j \le n-1} \alpha_{in}\alpha_{jn}}{\sum_{i=1}^{n-1} \alpha_{in}},
\]
the eigenvalue bound is immediate from, e.g., \cite[Thm~4.3.1]{MR1084815}.~\endproof

\section{A representation-theoretic view on Conjecture~\ref{conj}}
\label{sec:representation}
We now relate Conjecture~\ref{conj} to the 
representation theory of the symmetric group. Background on this theory is given 
in Appendix~\ref{sec:representationtheory}.
This section only contains standard results from representation theory, 
with a focus on transpositions of the symmetric group; 
see \cite[Section~3]{cesi} for a recent account in the context of the present paper.

We write $\rho^\lambda$ for Young's orthonormal irreducible representation corresponding to the partition $\lambda$, and set $V^\lambda_{ij}=I-\rho^\lambda_{ij}$.
Given the edge weights $\alpha =\{\alpha_{ij}\}$ of a graph $G$ with $n$ vertices, we set for $\lambda\vdash n$,
\[
\mathcal L^\lambda(\alpha) = \sum_{1\le i<j\le n} \alpha_{ij} V^\lambda_{ij}.
\]
Note that $\mathcal L^\lambda(\alpha)$ is a symmetric matrix, so it has real eigenvalues.
Write $\mu^\lambda(\alpha)$ for the eigenvalues of $\mathcal L^\lambda(\alpha)$, ordered so that $\mu_1^\lambda(\alpha)\le \ldots\le\mu^\lambda_{f^\lambda}(\alpha)$ where $f^\lambda$ is the dimension of the 
$V_{ij}^\lambda$:
the number of standard Young tableau with shape $\lambda$.

In the context of Markov chains arising from weighted graphs,
represention theory allows us to write their Laplacians as direct sums of matrices (up to a change of basis).
We first work this out for the weighted random walk, in which case the 
Laplacian is closely related to the so-called defining representation.

\begin{proposition}
\label{prop:laplaciandecomp}
There exists an $n\times n$ orthonormal matrix $S$ such that
for all weights $\alpha_{ij}$ on the edges of a graph with $n$ vertices,
\[
S \mathcal L_n^\RW(\alpha) S^\T = \mathcal L^{(n)}(\alpha) \oplus \mathcal L^{(n-1,1)}(\alpha).
\]
\end{proposition}

We note that the above decomposition of the random-walk Laplacian has a special structure.
Indeed, $\mathcal L^{(n)}(\alpha)$ equals zero regardless of the weights $\alpha$. Moreover,
$\mathcal L^{(n-1,1)}(\alpha)$ is closely related to $A_{n-1}$ reflection groups,
since $\rho^{(n-1,1)}_{ij}$ is the Householder reflection matrix 
corresponding to the transposition $(ij)$ of $\mathcal S_n$.
That is, each $\rho^{(n-1,1)}_{ij}$ acts on a point in $\R^{n-1}$ by reflecting it in a certain hyperplane.

\medskip
The Laplacian $\mathcal L_n^I(\alpha)$ of the interchange process,
which is an $n!\times n!$ matrix, can similarly be written as a direct sum (see also \cite[Section~3E]{MR964069}). 
This Laplacian is defined through
\[
\mathcal L_n^I(\alpha)(\sigma,\sigma') = 
\begin{cases}
-\alpha_{ij} & \text{if $\sigma' = (ij) \sigma$} \\
\sum_{1\le \ell< k \le n} \alpha_{\ell k} & \text{if $\sigma = \sigma'$} \\
0& \text{otherwise,}
\end{cases}
\]
where $\sigma,\sigma'\in\mathcal S_n$.
This Laplacian is closely related to the so-called regular representation.
We write $f^\lambda \mathcal L^\lambda(\alpha)$ for the direct sum 
of $f^\lambda$ copies of $\mathcal L^\lambda(\alpha)$.

\begin{proposition}
\label{prop:interchangelaplaciandecomp}
There exists an $n!\times n!$ orthonormal matrix $T$ such that
for all weights $\alpha_{ij}$ on the edges of a graph with $n$ vertices,
\[
T \mathcal L_n^I(\alpha)T^\T =
\bigoplus_{\lambda\vdash n} f^\lambda \mathcal L^\lambda(\alpha).
\]
\end{proposition}


The preceding proposition holds without any assumption on the sign of the weights $\alpha_{ij}$.
After defining signed weights on a graph with $k$ nodes through
\[
\alpha_{ij} = \begin{cases}
\gamma_i & \text{if $j=k$}\\
-\frac{\gamma_i\gamma_j}{\gamma_{1}+\ldots+\gamma_{k-1}}& \text{if $1\le i<j\le k-1$},
\end{cases}
\]
we immediately obtain a reformulation of Conjecture~\ref{conj} 
from Proposition~\ref{prop:interchangelaplaciandecomp}.

\begin{lemma}
\label{lem:alphaprime}
The following is equivalent to Conjecture~\ref{conj}.

Given any $k\ge 2$, 
the following holds for any $\lambda\vdash k$ and 
any nonnegative $\gamma_{1},\ldots,\gamma_{k-1}$:
\begin{equation}
\label{eq:conjrepr}
\sum_{i=1}^{k-1} \gamma_{i} V^\lambda_{ik} \ge 
\sum_{1\le i<j\le k-1} \frac{\gamma_{i}\gamma_{j}}{\gamma_{1}+\ldots+\gamma_{k-1}} V^\lambda_{ij}.
\end{equation}
\end{lemma}

\section{Conjecture~\ref{conj} implies Aldous' conjecture}
\label{sec:conjstronger}
In this section, we prove that Conjecture~\ref{conj} implies Aldous' conjecture.
We use mathematical induction on the number of vertices $n$. The conjecture trivially holds
if $n=2$.
Suppose Aldous' conjecture holds for all graphs with $n-1$ vertices. 


Consider an arbitrary weighted graph $G$ with $n$ vertices, 
and write $\alpha_{ij}\ge 0$ for the weight of edge $(i,j)$ in $G$.
By Proposition~\ref{prop:laplaciandecomp} and the fact that $\mu_1^{(n)}(\alpha)=0$,
Aldous' conjecture is that the second smallest eigenvalue of $\mathcal L_n^I(\alpha)$
equals $\mu_1^{(n-1,1)}(\alpha)$.
In view of Proposition~\ref{prop:interchangelaplaciandecomp}, 
this is equivalent with
\begin{equation}
\label{eq:toshow}
\mu^\lambda_1(\alpha) \ge \mu^{(n-1,1)}_1(\alpha )
\end{equation}
for all partitions $\lambda \vdash n$ with $\lambda\neq (n)$. This inequality trivially holds if $\lambda=(n-1,1)$, so we exclude this partition from further consideration.
Note that we do not need to assume that the graph $G$ be connected; since the right-hand side of
(\ref{eq:toshow}) vanishes if this is not the case, (\ref{eq:toshow}) then holds trivially since 
the $V_{ij}$ are positive semidefinite.

As before, we write $\alpha'$ for the weights given by (\ref{eq:defalphaprime}).
The induction hypothesis yields 
\[
\mu^{(n-2,1)}_1( \alpha' ) = \min_{\lambda' \,\vdash\, n-1: \,\lambda'\neq (n-1)} \mu^{\lambda'}_1(\alpha').
\]
To prove (\ref{eq:toshow}), we will show that the following string of inequalities holds if $\lambda\neq (n),(n-1,1)$: 
\begin{equation}
\label{eq:stringofineq}
\mu^{(n-1,1)}_1(\alpha) \le \mu^{(n-2,1)}_1(\alpha' ) = \min_{\lambda' \,\vdash\, n-1: \,\lambda'\neq (n-1)}  \mu^{\lambda'}_1(\alpha')  \le \mu^\lambda_1(\alpha).
\end{equation}
It is in the last inequality that we use Conjecture~\ref{conj}, but we first prove the first inequality.

\begin{lemma} 
\label{lem:firstineq}
We have
\[
\mu^{(n-1,1)}_1(\alpha) \le \mu^{(n-2,1)}_1(\alpha') \le \mu^{(n-1,1)}_2(\alpha) \le \ldots\le 
\mu^{(n-2,1)}_{n-2}(\alpha') \le \mu^{(n-1,1)}_{n-1}(\alpha).
\]
\end{lemma}
\proof
Consider the decomposition in Proposition~\ref{prop:laplaciandecomp}.
Since $\mathcal L_n^\RW(\alpha)$ and $\mathcal L_n^\RW(\alpha')$ are positive semidefinite
and $\mathcal L^{(n)}(\alpha)=\mathcal L^{(n)}(\alpha')$ are one-dimensional and
equal to zero, we conclude that $\mu_1^n(\alpha)=\mu_1^n(\alpha') =0$ and that
the largest $n-1$ eigenvalues of $\mathcal L_n^\RW(\alpha)$ and $\mathcal L_n^\RW(\alpha')$ 
are given by the eigenvalues of $\mathcal L^{(n-1,1)}(\alpha)$ and $\mathcal L^{(n-1,1)}(\alpha')$, 
respectively.
In particular, $\mu_{k+1}^n(\alpha) = \mu_k^{(n-1,1)}(\alpha)$ and 
$\mu_{k+1}^n(\alpha') = \mu_k^{(n-1,1)}(\alpha')$, and Proposition~\ref{prop:interlacing} yields
\[
\mu_1^{(n-1,1)}(\alpha')\le \mu_1^{(n-1,1)}(\alpha)\le \mu_2^{(n-1,1)}(\alpha')\le \ldots \le \mu_{n-1}^{(n-1,1)}(\alpha')\le \mu_{n-1}^{(n-1,1)}(\alpha).
\]
The claim readily follows from these inequalities, for instance
after noting that $\mu^{(n-1,1)}_1(\alpha')=\mu^{(n-1)}_1(\alpha')=0$ and 
$\mu^{(n-1,1)}_k(\alpha')=\mu^{(n-2,1)}_{k-1}(\alpha')$ for $k\ge 1$ 
by the branching rule (see Appendix~A).\endproof

\begin{lemma}
\label{lem:lastineq}
If Conjecture~\ref{conj} holds, then we have for $\lambda\ne (n),(n-1,1)$,
\[
\min_{\lambda' \,\vdash\, n-1: \,\lambda'\neq (n-1)}  \mu^{\lambda'}_1(\alpha')  \le \mu^\lambda_1(\alpha).
\]
\end{lemma}
\proof
Under Conjecture~\ref{conj}, we have by Lemma~\ref{lem:alphaprime},
\begin{eqnarray}
\mathcal L^\lambda(\alpha) &=& 
\sum_{1\le i<j\le n} \alpha_{ij} V_{ij}^\lambda =
\sum_{1\le i<j\le n-1} \alpha_{ij} V_{ij}^\lambda + \sum_{i=1}^{n-1} \alpha_{in} V_{in}^\lambda\nonumber\\ 
&\ge& \sum_{1\le i<j\le n-1} \alpha_{ij} V_{ij}^\lambda + 
\sum_{1\le i<j\le n-1} \frac{\alpha_{in}\alpha_{jn}}{\alpha_{1n}+\ldots+\alpha_{n-1,n}} V_{ij}^\lambda\label{eq:conjectureineq}\\
&=&\sum_{1\le i<j\le n-1} \alpha_{ij}' V_{ij}^\lambda = \mathcal L^\lambda(\alpha'),\nonumber
\end{eqnarray}
so that $\mu^\lambda_k(\alpha)\ge \mu^\lambda_k(\alpha')$ for all $k$.
By the branching rule (see Appendix~\ref{sec:representationtheory}), since $\alpha'_{in}=0$ for all $i$,
\[
\mathcal L^\lambda(\alpha') = \bigoplus_{\lambda'\vdash n-1 : \lambda' \nearrow \lambda} \mathcal L^{\lambda'}(\alpha').
\]
We conclude that
\[
\mu^\lambda_1(\alpha)\ge 
\mu^\lambda_1(\alpha') = \min_{\lambda'\vdash n-1 : \lambda' \nearrow \lambda} 
\mu^{\lambda'}_1(\alpha')\ge \min_{\lambda'\vdash n-1 : \lambda' \ne (n-1)} 
\mu^{\lambda'}_1(\alpha'),
\]
where the last inequality follows from $\lambda\neq (n),(n-1,1)$.
\endproof

\medskip
We have thus finished the proof that Conjecture~\ref{conj} implies Aldous' conjecture.
Before continuing, we mention a corollary of the proof which we use in 
Section~\ref{sec:application}.

\begin{corollary}
\label{cor:induction}
Suppose Conjecture~\ref{conj} holds for $k\le K$.
Let $\alpha_{ij}\ge 0$ be the weight of edge $(i,j)$
in a given weighted graph $G$ with $n$ vertices,
and suppose that 
at most $K-1$ of the $n-1$ possible $\alpha_{in}$ are strictly positive.

If Aldous' conjecture holds for the graph on $n-1$ vertices induced by edge weights $\alpha'$ given in
(\ref{eq:defalphaprime}),
then it holds for $G$.
\end{corollary}
\proof
The equality in (\ref{eq:stringofineq}) holds by assumption.
The first inequality in (\ref{eq:stringofineq}) holds by Lemma~\ref{lem:firstineq}, so we focus on the
last inequality and the proof of Lemma~\ref{lem:lastineq} in particular.

When at most $K-1$ of the $n-1$ possible $\alpha_{in}$ are strictly positive,
we may assume without loss of generality that these are $\alpha_{1n},\ldots,\alpha_{K-1,n}$.
By repeated application of the branching rule and a change of basis on interchanging $n$ and $K$,
we get
\[
\sum_{i=1}^{n-1} \alpha_{in} V_{in}^\lambda = \bigoplus_{\lambda' \vdash K: \lambda'\nearrow\ldots\nearrow\lambda} 
\left( \sum_{i=1}^{K-1} \alpha_{in} V_{in}^{\lambda'}\right)
\cong \bigoplus_{\lambda' \vdash K: \lambda'\nearrow\ldots\nearrow\lambda} 
\left( \sum_{i=1}^{K-1} \alpha_{in} V_{iK}^{\lambda'}\right),
\]
where the direct product should be taken over all possible simple paths from $\lambda'$ to $\lambda$
in the Hasse diagram of Young's lattice.
Thus, one can deduce (\ref{eq:conjectureineq}) from Conjecture~\ref{conj} with $k=K$
and the last inequality in (\ref{eq:stringofineq}) holds in that case.
\endproof

\medskip
This corollary is of particular interest when $K=2$, in which case Conjecture~\ref{conj} holds trivially.
If the $n$-th vertex is incident to exactly one other vertex, then
the modified weights $\alpha'$ on the edges of the `small' graph
consisting of the vertices $1,\ldots,n-1$ are simply equal to the unmodified weights $\alpha$.
In this context, Corollary~\ref{cor:induction} is essentially equivalent to the induction
step in Handjani and Jungreis~\cite{MR1419872}, and it readily implies that Aldous' conjecture holds for trees.
A related argument is given by Cesi~\cite[Sec.~3]{cesi2}.

\section{Proof of Conjecture~\ref{conj} for $k\le 4$}
\label{sec:proofconjS4}
Conjecture~\ref{conj} for $k=4$ implies the conjecture for $k<4$, so this section focuses on proving 
Conjecture~\ref{conj} for $k=4$.
In view of Lemma~\ref{lem:alphaprime}, it suffices
to prove (\ref{eq:conjrepr}) for all $\lambda\vdash 4$.
We do so by making use of the explicit forms of Young's orthonormal
irreducible representations given in Appendix~\ref{sec:reprS4}.
In particular, we use the vectors $v^\lambda_{ij}$ introduced in Appendix~\ref{sec:reprS4}.
Throughout this section, for notational convenience, we suppress
the superscripts $\lambda$ in these vectors, so that,
e.g., $v_{ij}$ stands for $v_{ij}^{(3,1)}$ in Section~\ref{sec:repr31} while 
it stands for $v_{ij}^{(2,2)}$ in Section~\ref{sec:repr22}. 
We follow the same notational convention for the superscripts $\lambda$ in $V^\lambda_{ij}$.

\subsection{$\lambda={(4)}$}
Since $V_{ij} = 0$ for $1\le i<j\le 4$, (\ref{eq:conjrepr}) trivially holds.

\subsection{$\lambda={(3,1)}$}
\label{sec:repr31}
Since $v_{ij} = v_{i4}-v_{j4}$ for $1\le i<j\le 3$, 
we readily find that
\[
(\gamma_1+\gamma_2+\gamma_3) \sum_{i=1}^3 \gamma_i V_{i4} - \sum_{1\le i<j\le 3}\gamma_i\gamma_j V_{ij}
= (\gamma_1 v_{14} + \gamma_2 v_{24}+ \gamma_3 v_{34}) (\gamma_1 v_{14} + \gamma_2 v_{24}+ \gamma_3 v_{34})^\T,
\]
as in the proof of Proposition~\ref{prop:interlacing}.
This implies (\ref{eq:conjrepr}) for $\lambda=(3,1)$.

\subsection{$\lambda=(2,2)$}
\label{sec:repr22}
We first show that  
\begin{eqnarray}
\lefteqn{(\gamma_1+\gamma_2+\gamma_3) \sum_{i=1}^3 \gamma_i V_{i4} - \sum_{1\le i<j\le 3}\gamma_i\gamma_j V_{ij}} \nonumber\\
&=& \sum_{1\le i<j\le 3} (\gamma_i v_{i4} - (-1)^{i-j} \gamma_j v_{j4}) (\gamma_i v_{i4} - (-1)^{i-j} \gamma_j v_{j4})^\T
- \sum_{i=1}^3 \gamma_i^2 V_{i4}.
\label{eq:toshow22}
\end{eqnarray}
Since $v_{12} = -(v_{14}-v_{24})$,
we see that
\begin{eqnarray*}
(\gamma_1+\gamma_2)(\gamma_1 V_{14} + \gamma_2 V_{24})-\gamma_1\gamma_2 V_{12} &=& 
(\gamma_1+\gamma_2)(\gamma_1 V_{14} + \gamma_2 V_{24})-\gamma_1\gamma_2 (v_{14}-v_{24})(v_{14}-v_{24})^\T \\
&=& (\gamma_1 v_{14} +\gamma_2 v_{24}) (\gamma_1 v_{14} +\gamma_2 v_{24})^\T.
\end{eqnarray*}
After two similar calculations using $v_{13} = v_{14}+v_{34}$ and $v_{23}=v_{24}-v_{34}$,
we find that for $1\le i<j\le 3$,
\[
(\gamma_i+\gamma_j)(\gamma_i V_{i4} + \gamma_j V_{j4})-\gamma_i\gamma_j V_{ij}=
(\gamma_i v_{i4} - (-1)^{i-j} \gamma_j v_{j4}) (\gamma_i v_{i4} - (-1)^{i-j} \gamma_j v_{j4})^\T.
\]
After summing these identities over $1\le i<j\le 3$ and some rearranging, we get (\ref{eq:toshow22}).

To show that the right-hand side of (\ref{eq:toshow22}) is positive semidefinite,
we first introduce
\[
u = 
\left(\begin{array}{c}
\gamma_1\gamma_2 \\ -\gamma_1\gamma_3 \\ \gamma_2\gamma_3 
\end{array}\right),
\quad 
w = 
\left(\begin{array}{c}
1 \\ -1 \\ 1 
\end{array}\right),
\]
and write $P_u=I_3-u u^\T/\|u\|^2$ for the projection matrix on the hyperplane orthogonal to $u$.
By the Courant-Fischer variational theorem,
the smallest eigenvalue of the matrix in (\ref{eq:toshow22}) equals
\begin{eqnarray*}
\lefteqn{\min_{x\in\R^2: \|x\|=1} 
\sum_{1\le i<j\le 3} [\gamma_i x^\T v_{i4} -(-1)^{i-j} \gamma_j x^\T v_{j4}]^2
- \sum_{i=1}^3 [\gamma_i x^\T v_{i4}]^2} \\
&\ge& \min_{z\in\R^3: \|z\| =1,  u^\T z=0} [(z_1+z_2)^2+(z_1-z_3)^2+(z_2+z_3)^2 -z_1^2-z_2^2-z_3^2] \\
&=& \min_{z\in\R^3: \|z\| =1,  u^\T z=0} z^\T\left[
2I_3 - w w^\T \right]z 
= 2-\max_{z\in\R^3: \|z\| =1,  u^\T z=0} (z^\T w)^2 \\
&=& 2-\max_{z\in\R^3: \|z\| =1} (z^\T P_u w)^2 
= 2-\tr (P_u w w^\T P_u)\\
&=& 2-w^\T P_u w
= -1 + (w^\T u)^2/\|u\|^2,
\end{eqnarray*}
where the first inequality follows from $v_{14}-v_{24}+v_{34} =0$.
Since $(w^\T u)^2 \ge \|u\|^2$, we have proven (\ref{eq:conjrepr}) for $\lambda={(2,2)}$.

\subsection{$\lambda={(2,1^2)}$}
Since $V_{ij} = 2I_3-v_{ij}v_{ij}^\T$ and 
$v_{ij} = v_{i4}-v_{j4}$ for $1\le i<j\le 3$, we find that
\begin{eqnarray*}
\lefteqn{(\gamma_1+\gamma_2+\gamma_3) \sum_{i=1}^3 \gamma_i V_{i4} - \sum_{1\le i<j\le 3}\gamma_i\gamma_j V_{ij}}\\
&=&2\left(\sum_{1\le i\le j\le 3} \gamma_i\gamma_j\right)I_3 - 
(\gamma_1+\gamma_2+\gamma_3) \sum_{i=1}^3 \gamma_i v_{i4}v_{i4}^\T+
\sum_{1\le i<j\le 3}\gamma_i\gamma_j v_{ij}v_{ij}^\T\\
&=&2\left(\sum_{1\le i\le j\le 3} \gamma_i\gamma_j\right)I_3 - 
(\gamma_1 v_{14} + \gamma_2 v_{24} + \gamma_3 v_{34})(\gamma_1 v_{14} + \gamma_2 v_{24} + \gamma_3 v_{34})^\T,
\end{eqnarray*}
where the last equality follows from the same calculation as in Section~\ref{sec:repr31}.
Therefore, the smallest eigenvalue of this matrix is
\[
2\left(\sum_{1\le i\le j\le 3} \gamma_i\gamma_j\right) - 
(\gamma_1 v_{14} + \gamma_2 v_{24} + \gamma_3 v_{34})^\T(\gamma_1 v_{14} + \gamma_2 v_{24} + \gamma_3 v_{34}),
\]
which equals zero since $v_{i4}^\T v_{j4}=1$ for $i\neq j$ while $v_{i4}^\T v_{i4}=2$.
This proves (\ref{eq:conjrepr}) for $\lambda=(2,1^2)$.

\subsection{$\lambda={(1^4)}$}
Since $V_{ij} = 2$ for all $i,j$, 
(\ref{eq:conjrepr}) reduces to $\sum_{i=1}^{k-1} \gamma_i^2 + \sum_{1\le i<j\le k-1} \gamma_i\gamma_j\ge 0$,
which clearly holds.

\section{Proof of Conjecture~\ref{conj} for $\gamma_{1}=\ldots=\gamma_{k-1}$}
\label{sec:proofconjJM}
This section proves Conjecture~\ref{conj} for $\gamma_{1}=\ldots=\gamma_{k-1}$.
The main ingredients are Jucys-Murphy matrices and a content minimization calculation
for standard Young tableaux, see
Appendix~\ref{sec:representationtheory} for definitions.

Fix $k\ge 2$ and choose some partition $\lambda\vdash k$. We need to prove that 
\begin{equation}
\label{eq:JMdiagmatrix}
(k-1) \sum_{i=1}^{k-1} V^\lambda_{ik} -
\sum_{1\le i<j\le k-1}  V^\lambda_{ij}\ge 0.
\end{equation}
The left-hand side of (\ref{eq:JMdiagmatrix}) can be written in terms of 
Jucys-Murphy matrices as
\[
\frac 12 k(k-1) I_{f^\lambda}
 + \sum_{i=1}^{k} X_{i}^\lambda-k X^\lambda_{k}.
\]
In particular, it is a diagonal matrix
and its diagonal elements are readily found.
Indeed, element $(t,t)$ of this matrix is calculated from tableau $t$ through
\[
\frac 12 k(k-1) + \sum_{i=1}^{k} c^t_i - k c^t_k,
\]
where $c^t_i$ is the content of the box containing $i$ in tableau $t$.
The sum over $i$ is the sum of all contents corresponding to a Young tableau of shape $\lambda$,
which can be expressed in terms of $\lambda$ by noting that the sum over the contents
in the $j$-th row equals $\frac 12 \lambda_j(\lambda_j-1) - (j-1)\lambda_j$.
As this is independent of $t$, the smallest 
diagonal element of the matrix on the left-hand side of (\ref{eq:JMdiagmatrix}) 
corresponds to a tableau $t$ for which $c^t_k$ is maximized, i.e., to a tableau $t$ with $c^t_k = \lambda_1-1$.
Therefore,
we find that the smallest eigenvalue of (\ref{eq:JMdiagmatrix}) equals
\begin{eqnarray*}
\lefteqn{\frac 12 k(k-1) + \sum_{i=1}^{k} c^t_i - k (\lambda_1-1)}\\
&=& 
\frac 12 k(k-1) + \sum_{j=1}^\infty \frac 12 \lambda_j(\lambda_j-1) 
 - \sum_{j=1}^\infty (j-1)\lambda_j - k (\lambda_1-1)\\
&=& \frac 12 k^2 + \sum_{j=1}^\infty \frac 12 \lambda_j^2 
- \sum_{j=1}^\infty (j-1)\lambda_j - k \lambda_1 = \frac 12 (k-\lambda_1)^2 +
\frac 12 \sum_{j=2}^\infty  \lambda_j^2 - \sum_{j=1}^\infty (j-1)\lambda_j
\\&=& \frac 12 \left(\sum_{j=2}^\infty \lambda_j\right)^2 +
\frac 12 \sum_{j=2}^\infty  \lambda_j^2 - \sum_{j=1}^\infty (j-1)\lambda_j
= \sum_{j=2}^\infty \lambda_j^2 + \sum_{j=2}^\infty \sum_{i=2}^{j-1}\lambda_i\lambda_j 
-\sum_{j=1}^\infty (j-1)\lambda_j\\
&\ge& \sum_{j=2}^\infty \lambda_j^2 + \sum_{j=2}^\infty (j-2)\lambda_j^2 
-\sum_{j=1}^\infty (j-1)\lambda_j
= \sum_{j=2}^\infty (j-1)\lambda_j (\lambda_j-1).
\end{eqnarray*}
Since each $\lambda_j$ is a nonnegative integer, this is clearly nonnegative.
This proves Conjecture~\ref{conj} for $\gamma_{1}=\ldots=\gamma_{k-1}$.

\section{Weighted graphs with nested triangulation}
\label{sec:application}
In this section, we introduce a class of weighted graphs for which we prove Aldous' conjecture.
This class includes all trees and all cycles of arbitrary length,
and it arises by repeated application of Corollary~\ref{cor:induction} for $K=4$.

\begin{figure}
\resizebox{115pt}{!}{\includegraphics*{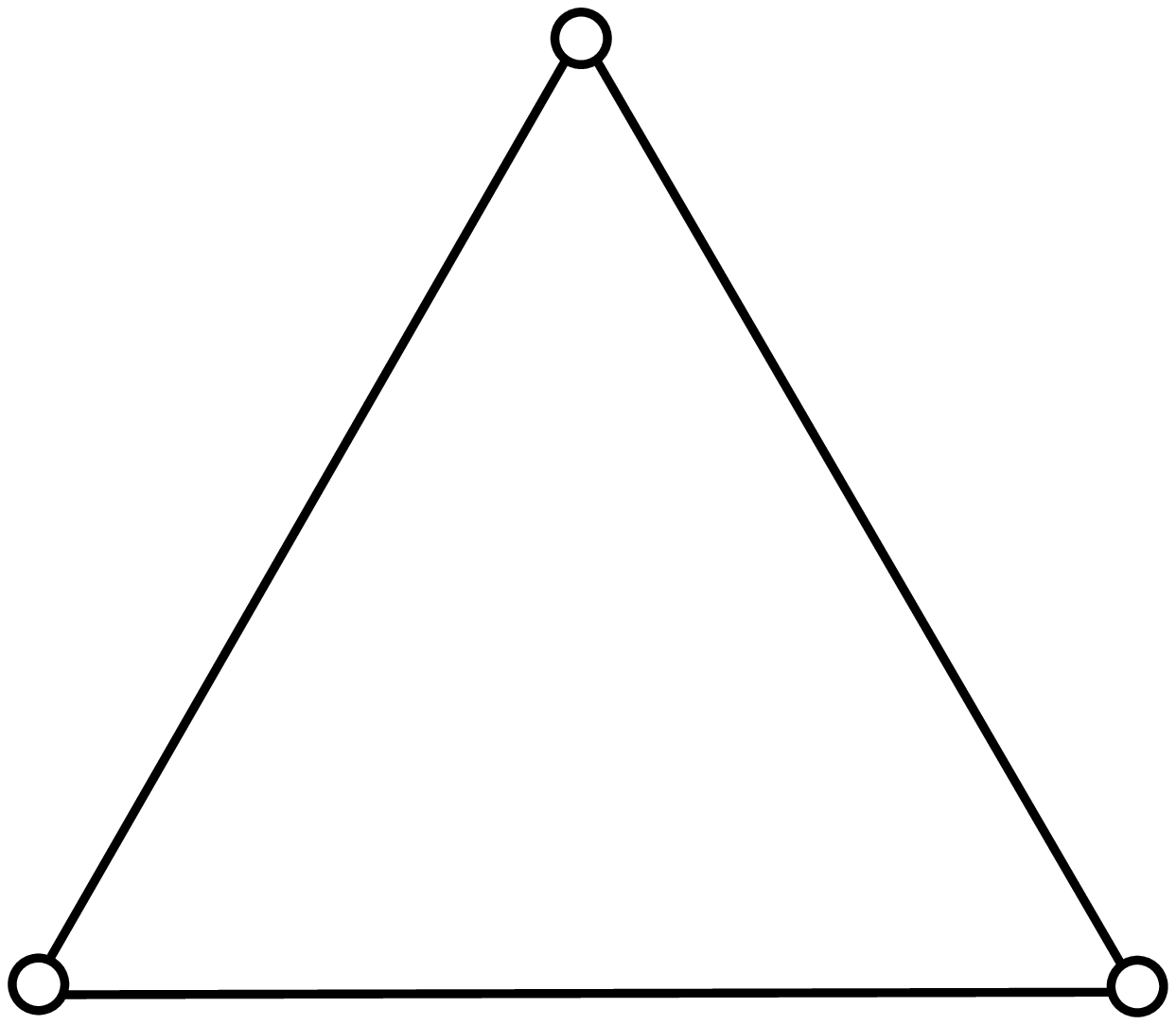}}
\hspace{7mm}
\resizebox{115pt}{!}{\includegraphics*{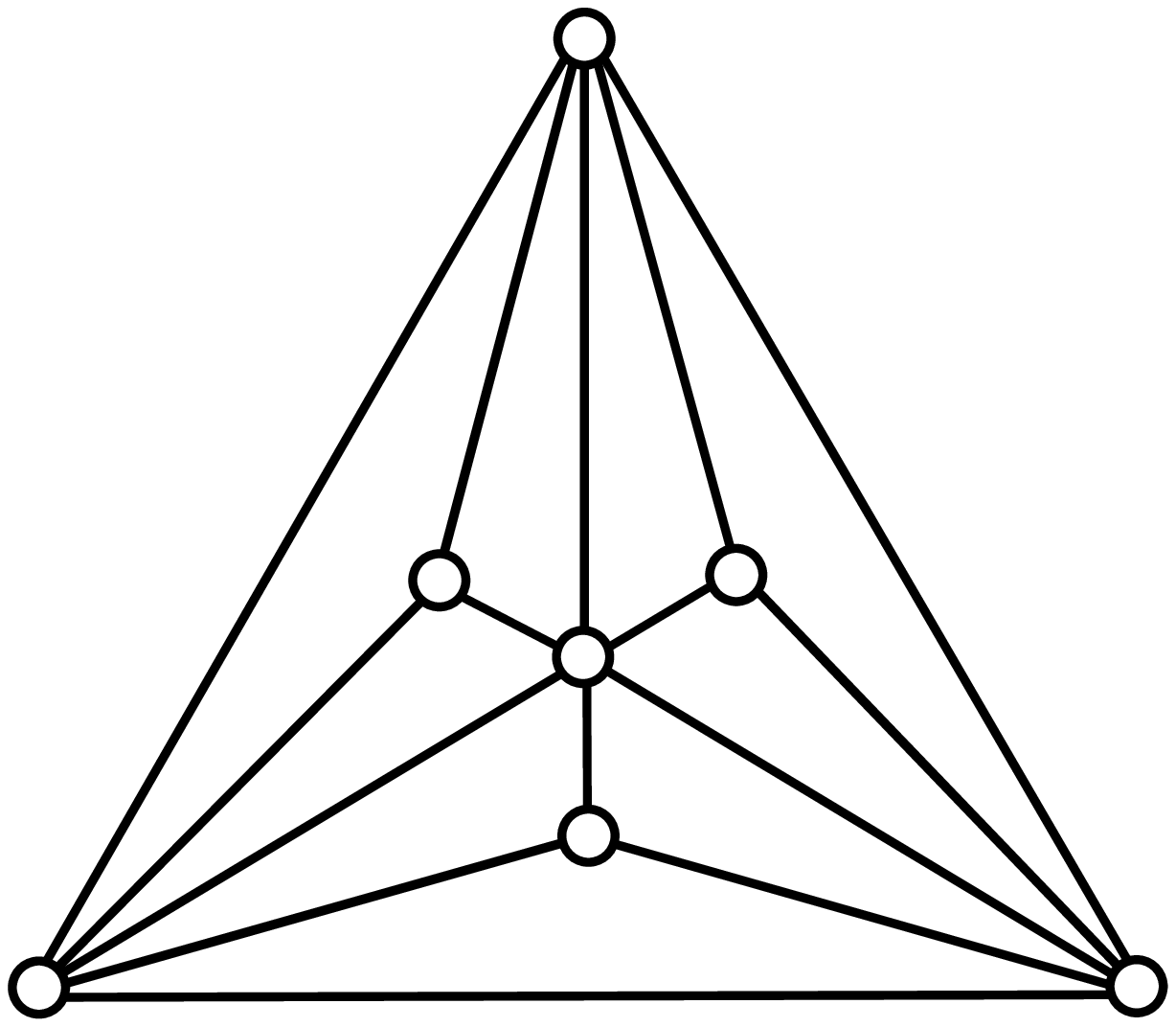}}
\hspace{7mm}
\resizebox{115pt}{!}{\includegraphics*{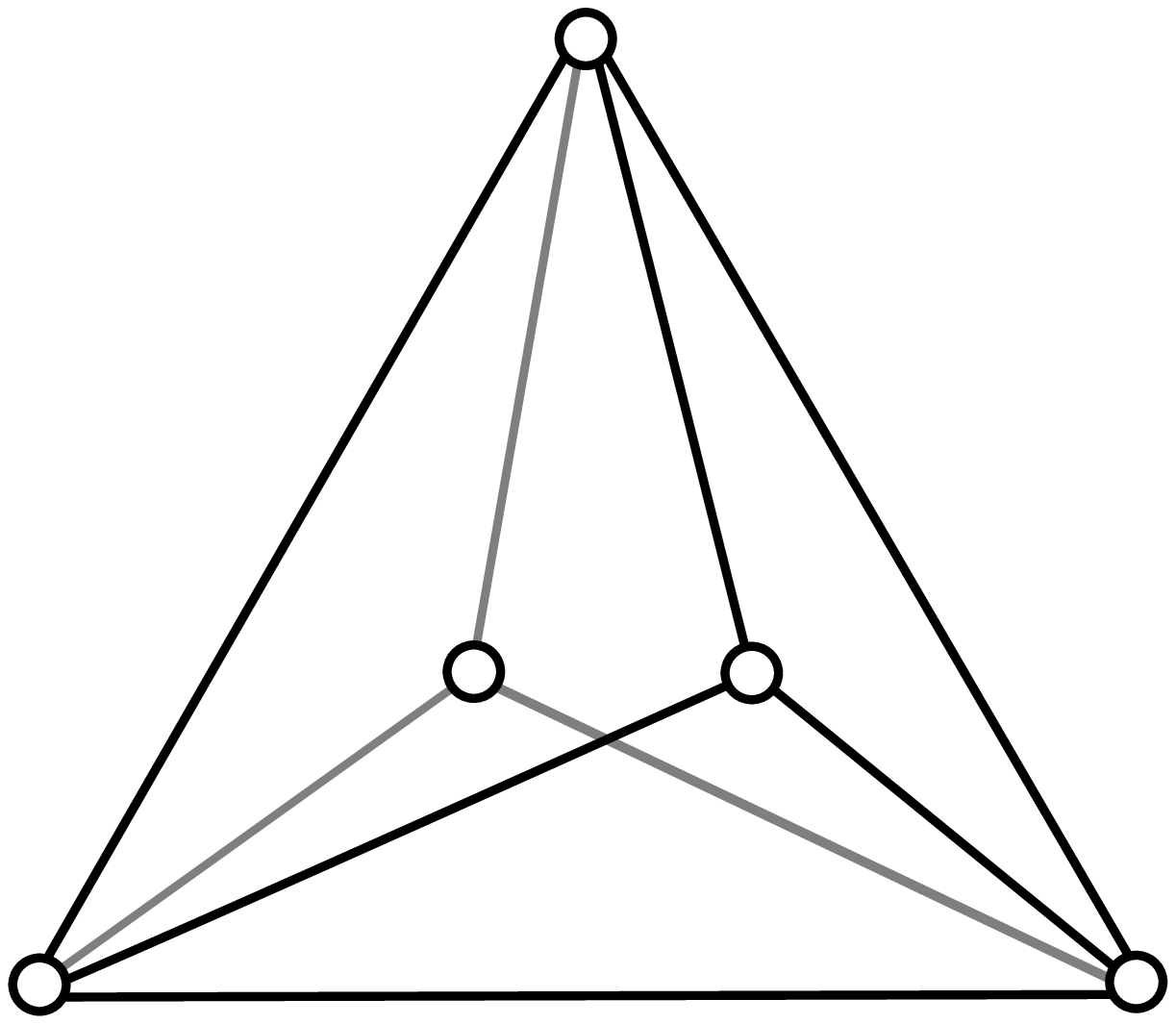}}
\caption{\label{fig:graphs} The graphs $T_{0,N}=K_3$ (left), $T_{2,1}$ (center), and $T_{1,2}$ (right).}
\end{figure}
Our graphs with nested triangulations  
are parameterized by two integers:
a branching parameter $N\ge 1$ and a depth parameter $D\ge 0$.
For a given $N$, the graphs $\{T_{i,N}: i\ge 0\}$ are nested in the sense that 
$T_{i,N}$ is a subgraph of $T_{i+1,N}$ for $i\ge 0$. 
The graphs are defined recursively as follows.
Let $T_{0,N}$ be the complete graph on 3 vertices: $T_{0,N}=K_3$.
For each cycle of length 3 that is present in $T_{i,N}$ but not in $T_{i-1,N}$,
we construct $T_{i+1,N}$ by adding
$N$ vertices to $T_{i,N}$, and by adding 3 new edges for each new vertex to connect it
to the 3 vertices of the given cycle. Thus, $3N$ edges are added for each cycle
of length 3 in $T_{i,N}$ but not in $T_{i-1,N}$.
The vertices of $T_{D,N}$
can be partitioned into $D+1$ levels according to the stage at which they have been
added. Examples are given in Figure~\ref{fig:graphs}. Note that $T_{D,1}$ is
a maximal planar (triangulated) graph for any $D\ge 1$, but that not all maximal planar graphs are 
graphs with nested triangulations.
Also note that $T_{3,1}$ has $K_{3,3}$, the complete bipartite graph on six vertices, as a subgraph and
it is therefore nonplanar by Kuratowski's theorem.



\begin{proposition}
\label{prop:TDN}
For any $D\ge 0$, $N\ge 1$, let $T_{D,N}$ have arbitrary nonnegative interchange rates on its edges
and assume that the graph remains connected after removing zero-rate edges.
Aldous' conjecture holds for this graph.
\end{proposition}
\proof
We use induction.
Since Aldous' conjecture trivially holds for a connected graph with two vertices, 
we conclude from Corollary~\ref{cor:induction} that Aldous' conjecture holds for the triangle $K_3$. 
Since each vertex at level $i+1$ is incident to exactly 3 vertices at lower levels, we may
repeatedly use Corollary~\ref{cor:induction} to deduce the claim for $T_{i+1,N}$ from
the claim for $T_{i,N}$.~\endproof

\medskip
This proposition is of particular interest for $D=N=1$, in which case $T_{D,N}$ is the complete graph $K_4$ on
four vertices. Proposition~\ref{prop:TDN} then states that Aldous' conjecture holds for {\em all}
weighted graphs with four vertices.

Choosing some of the interchange rates equal to zero in Proposition~\ref{prop:TDN} proves
Aldous' conjecture for some special classes of graphs.
For instance, any tree can be embedded in a graph with nested triangulations.
Indeed, given any tree, let $D$ be the maximum distance to the root and let $N$ be the maximum degree.
It is readily seen that one can embed the tree into $T_{D,N}$ by mapping a vertex at distance $i\ge 0$ 
from the root to a vertex at level $i$ in $T_{D,N}$.
Thus, Proposition~\ref{prop:TDN} recovers the main result from \cite{MR1419872} in this case.

Instead of showing that a given graph is a subgraph of a graph with inner triangulations 
and appealing to Proposition~\ref{prop:TDN},
the following prodecure is an alternative for showing that Aldous' conjecture must hold
according to the results of this paper.
Corollary~\ref{cor:induction} implies that Aldous' conjecture holds for graphs 
which (after removing all edge weights) can be reduced to an edge by repeatedly using the 
following permissible rules:
\begin{itemize}
\item Degree-one reduction: delete a degree-one vertex and its incident edge.
\item Series reduction: delete a degree-two vertex $k$ and its two incident edges $(i,k)$ and $(j,k)$, and add in a new edge $(i,j)$.
\item Parallel reduction: delete one of a pair of parallel edges.
\item Y-$\Delta$ transformation: delete a vertex $k$ and its three incident
edges $(i,k)$, $(j,k)$, $(\ell, k)$ and add in a triangle $ij\ell$.
\end{itemize}
These operations also appear in the context of 
star-triangle reducibility of a graph~\cite{MR1222626}, but it is important to note that the $\Delta$-Y transformation 
(which is the inverse of the Y-$\Delta$ transformation)
is not permissible here.

\begin{figure}
\resizebox{150pt}{!}{\includegraphics*{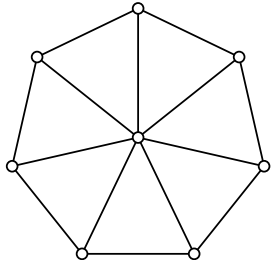}}
\caption{\label{fig:wheel} The wheel graph $W_7$ with 7 vertices.}
\end{figure}

Wheel graphs are examples of graphs which can be reduced to an edge 
using these operations. We write $W_n$ for the wheel graph with $n$ vertices, see 
Figure~\ref{fig:wheel} for $W_7$.
Indeed, one obtains $W_{n-1}$ from $W_n$ after applying a Y-$\Delta$ transformation to 
one of the outer vertices of $W_n$ followed by three parallel reductions. This procedure 
can be repeated until $W_3=K_3$ arises, which is readily reduced to an edge.
Note that cycles are subgraphs of wheel graphs:
choose the interchange rates on the spokes of the wheel equal to zero except for two adjacent spokes, 
and also let the interchange rate vanish on the edge incident to the two outer vertices of these two spokes.
Thus, we have also proven that Aldous' conjecture holds for weighted cycles.

\section{Discussion}
\subsection*{Other Markov processes with the same spectral gap as the random walk} 
Apart from the interchange process and the random walk,
several other natural Markov chains arise from the interchange dynamics on a weighted graph. 
Indeed, one may allow several vertices to receive the same label, which 
can be thought of as a color.
Interchanging nodes with the same color then does not change the color configuration on the 
graph. Thus, for each possible initial configuration of colors,
one obtains a continuous-time Markov chain.
One can think of these processes as parameterized by Young diagrams (partitions), where 
each row corresponds to a color and the number of boxes in each row correspond to 
the number of vertices to receive this color.
The resulting process can be interpreted as a random walk on a so-called Schreier graph,
see also Cesi~\cite{cesi2}.
The interchange process is a special case of this construction
with $\lambda=(1^n)$, i.e., all vertices have different colors.
Similarly, the random walk process arises on setting $\lambda=(n-1,1)$, i.e., one vertex has a different color from the other vertices.

After a change of basis, the intensity matrices for these Markov processes can be 
written as a direct sum of irreducible representations as in
Propositions~\ref{prop:laplaciandecomp} and \ref{prop:interchangelaplaciandecomp}.
This is called Young's rule;
indeed, the intensity matrix corresponding to partition $\lambda$ naturally arises 
from the $M^\lambda$ module in representation theory.
The multiplicities of the irreducible representations are given by the so-called Kostka numbers.
As a consequence of the resulting block structure, all of the intensity matrices 
(except for the trivial one corresponding to $\lambda=(n)$)
contain the irreducible representation corresponding to the partition $(n-1,1)$.
Thus, if Conjecture~\ref{conj} can be shown to hold, all of these processes 
have exactly the same spectral gap.

\subsection*{Gelfand-Tsetlin patterns}
By Proposition~\ref{prop:interlacing}, subsequent removal of vertices 
and updating of the weights according to (\ref{eq:defalphaprime}) yields a 
Gelfand-Tsetlin pattern, i.e., a collection of subsequent interlaced sequences.
Subsequent removal of vertices {\em without} weight updating yields a nondecreasing spectral-gap sequence, 
an observation which has previously proven useful in the context of 
Aldous' conjecture \cite{MR2415139,starrconomos}.
The significance of the Gelfand-Tsetlin structure
is currently unclear.

\subsection*{The cut-off phenomenon}
It is a natural question whether the results of this paper can be exploited
to study the cut-off phenomenon for Markov chains.
This question is currently open.
Proving a cut-off phenomenon requires control over the whole spectrum, not only near the edge.  
A variety of known results \cite{MR2023654}, e.g., on $\ell$-adjacent transposition walks, 
suggests that (pre)cut-off thresholds for interchange
processes have an extra $\log(n)$ factor when compared to the corresponding random walk processes.
Proposition~\ref{prop:interchangelaplaciandecomp} suggests 
that the second smallest eigenvalue of the Laplacian of the interchange process
typically has multiplicity $f^{(n-1,1)} = n-1$ if Aldous' conjecture holds,
which may explain the additional $\log(n)$ factor.

\subsection*{Electric networks}
Section~\ref{sec:application} showed how the Y-$\Delta$ transformation naturally arises in
the context of Corollary~\ref{cor:induction}, but there may be a deeper connection.
For $n=4$, the definition of $\alpha'$ in (\ref{eq:defalphaprime}) in terms of $\alpha$
appears in formulas for the resistance in electrical networks when a $\Delta$ is transformed
into a Y.
The recent work of Caputo {\em et al.}~\cite{caputo} sheds some light on this.

\section*{Acknowledgments}
The author would like to thank Prasad Tetali and David Goldberg for valuable discussions,
and Kavita Ramanan for helpful comments on an earlier draft.

\appendix
\section{Background on representation theory}
\label{sec:representationtheory}
This section reviews the elements of representation theory used in the body of this paper.
More comprehensive accounts can be found in \cite{MR1824028,fulton:youngtableaux1997,MR513828,MR1443185}.

A {\em partition} $\lambda$ of $n$, written $\lambda\vdash n$, is a sequence of 
nonnegative integers $(\lambda_1,\lambda_2,\ldots)$ with $\sum_{j=1}^\infty \lambda_j=n$ 
and $\lambda_1\ge\lambda_2\ge \ldots$.
For notational convenience, we suppress the zero elements of the sequence.
Also, if the integer $k$ appears $m$ times in the partition $\lambda$, we replace the $m$ copies of
$k$ by a single copy of $k^m$.
For instance, $(4,3^2,1)$ is shorthand for $(4,3,3,1,0,\ldots)$.
A partition can be identified with a Young diagram, which is 
a collection of $n$ boxes arranged in left-justified rows, with the $i$-th row
containing $\lambda_i$ boxes. 
For instance,
\[
(4,2,1) = \yng(4,2,1).
\]
Given two partitions $\lambda'$ and $\lambda$,
we write $\lambda' \nearrow \lambda$ if the Young diagram of $\lambda$
can be obtained from the Young diagram of $\lambda'$ by adding a box. 
For instance, we have
\[
\yng(4,2,1) \quad\nearrow \quad\yng (4,3,1),
\]
since a box is added to the second row. 
A natural related partial order on partitions is defined by diagram containment.
The set of all partitions equipped with this partial order is called {\em Young's lattice}.

Given a partition $\lambda$, a {\em standard Young tableau} with shape $\lambda\vdash n$ 
is the Young diagram corresponding to $\lambda$ with each of the numbers $1,\ldots,n$
inside one of the $n$ boxes,
in such a way that the numbers in each of the rows as well as in each of the columns of
the Young diagram are increasing. For instance,
\[
\young(1234,56,7) \quad {\rm and} \quad
\young(1357,24,6)
\]
are both standard Young tableaux with shape $(4,2,1)$. 
We write $f^\lambda$ for the number of different standard Young tableau with shape $\lambda$.
Note that a Young tableau with shape $\lambda$ can alternatively be thought of
as a saturated chain in Young's lattice starting from the empty partition \cite[Prop.~7.10.3]{stanley:ec2}.

The {\em content} of a box in a Young diagram is defined as the $x$-coordinate of 
the box minus its $y$-coordinate. 
Thus, the boxes in the following Young diagram contain their content:
\[
\young(0123,{\negone}0,{\negtwo}).
\]
For instance, in the two given examples of standard Young tableaux,
the box containing $7$ has content $-2$ and $3$, respectively.
We write $c^t_i$ for the content of the box containing $i$ in tableau $t$.

We next introduce Young's orthonormal irreducible representation corresponding to the partition $\lambda\vdash n$.
This is a family of $f^\lambda\times f^\lambda$ matrices 
$\{\rho^\lambda_\sigma: \sigma\in\mathcal S_n\}$ parameterized by elements of $\mathcal S_n$,
such that the matrices behave in the same way as the elements of $\mathcal S_n$ when multiplied together
(i.e., the mapping $\sigma\mapsto \rho^\lambda_\sigma$ is a group homomorphism).
When $\sigma$ is a transposition, say $\sigma=(ij)$, we 
write $\rho^\lambda_{ij}$ instead of $\rho^\lambda_{(ij)}$.
We let the standard basis vectors correspond to the standard Young tableaux with shape $\lambda$
under the dictionary (total) order, i.e., if
the numbers are read from left to right by rows, starting at the top row, the
first digit in which two tableaux disagree will be larger for the larger tableau.
We fix $2\le i\le n$ and specify $\rho^\lambda_{i-1,i}$; since the adjacent transpositions 
generate $\mathcal S_n$, this specifies the whole group representation:
\begin{itemize}
\item If $i$ and $i + 1$ are in the same row of $t$, then $\rho^\lambda_{i,i+1}(t,t)=1$.
\item If $i$ and $i + 1$ are in the same column of $t$, then $\rho^\lambda_{i,i+1} (t,t)=-1$.
\item Suppose $i$ and $i + 1$ are not in the same row or column of $t$.
Write $s$ for the standard Young tableau resulting from swapping $i$ and $i + 1$ in $t$.
Then we have
\[
\left(\begin{array}{cc}
\rho^\lambda_{i,i+1}(t,t) & \rho^\lambda_{i,i+1}(t,s) \\
\rho^\lambda_{i,i+1}(s,t) & \rho^\lambda_{i,i+1}(s,s)
\end{array}\right) =
\left(\begin{array}{cc}
r^{-1} & \sqrt{1-r^{-2}}\\
\sqrt{1-r^{-2}} & - r^{-1}
\end{array}\right),
\]
where $r = c^t_{i+1}-c^t_i$, the {\em axial distance} between the boxes containing $i$ and $i+1$. 
\end{itemize}
The elements of $\rho^\lambda_{i,i+1}$ left unspecified by these three rules are zero.
A {\em branching rule} holds, which reduces to
\[
\rho_{ij}^\lambda = \bigoplus_{\lambda'\vdash n-1 : \lambda' \nearrow \lambda} \rho_{ij}^{\lambda'}
\]
for transpositions $(ij)$ with $i<j<n$.

Of special importance in representation theory are the so-called {\em Jucys-Murphy elements};
they play a key role in the Vershik-Okounov approach to representation theory \cite{MR1443185}.
For $2\le j\le f^\lambda$, the Jucys-Murphy matrix corresponding to the partition $\lambda$ is defined as
$X^\lambda_j = \sum_{i=1}^{j-1} \rho^\lambda_{ij}$.
Their significance stems from the fact that these matrices commute; in fact, they are diagonal matrices.
Element $(t,t)$ of $X^\lambda_j$ equals $c_j^t$, the content of the box containing $j$ in tableau $t$.

\section{Young's orthonormal irreducible representations of $\mathcal S_4$}
\label{sec:reprS4}
This appendix evaluates Young's orthonormal irreducible representations of $\mathcal S_4$ at transpositions,
which is a key ingredient in Section~\ref{sec:proofconjS4}.
The given formulas can be verified with the definition of Young's orthonormal irreducible representation
in Appendix~\ref{sec:representationtheory}.

\noindent
For $1\le i<j\le 4$, we have $\rho^{(4)}_{ij}= 1$ (one-dimensional).

\noindent
For $1\le i<j\le 4$, we have
\[
\rho^{(3,1)}_{ij} = I_3 - v^{(3,1)}_{ij} {v^{(3,1)}_{ij}}^\T,
\]
where
\begin{eqnarray*}
&&v^{(3,1)}_{12} = \left(\begin{array}{c} 0 \\ 0\\\sqrt{2} \end{array}\right),
v^{(3,1)}_{13} = \left(\begin{array}{c} 0 \\ \sqrt{3/2}\\ \sqrt{1/2} \end{array}\right),
v^{(3,1)}_{14} = \left(\begin{array}{c} \sqrt{4/3} \\ \sqrt{1/6} \\ \sqrt{1/2} \end{array}\right), \\
&&v^{(3,1)}_{23} = \left(\begin{array}{c} 0 \\ \sqrt{3/2}\\ -\sqrt{1/2} \end{array}\right),
v^{(3,1)}_{24} = \left(\begin{array}{c} \sqrt{4/3} \\ \sqrt{1/6} \\ -\sqrt{1/2} \end{array}\right),
v^{(3,1)}_{34} = \left(\begin{array}{c} \sqrt{4/3} \\ -\sqrt{2/3}\\0 \end{array}\right).
\end{eqnarray*}

\noindent
For $1\le i<j\le 4$, we have
\[
\rho^{(2,2)}_{ij} = I_2 - v^{(2,2)}_{ij} {v^{(2,2)}_{ij}}^\T,
\]
where
\begin{eqnarray*}
&&v^{(2,2)}_{12} = \left(\begin{array}{c}  0\\\sqrt{2} \end{array}\right),
v^{(2,2)}_{13} = \left(\begin{array}{c}  \sqrt{3/2}\\ \sqrt{1/2} \end{array}\right),
v^{(2,2)}_{14} = \left(\begin{array}{c} \sqrt{3/2}\\ -\sqrt{1/2} \end{array}\right), \\
&&v^{(2,2)}_{23} = \left(\begin{array}{c} \sqrt{3/2}\\ -\sqrt{1/2} \end{array}\right),
v^{(2,2)}_{24} = \left(\begin{array}{c} \sqrt{3/2} \\ \sqrt{1/2} \end{array}\right),
v^{(2,2)}_{34} = \left(\begin{array}{c} 0\\ \sqrt{2} \end{array}\right).
\end{eqnarray*}

\noindent
For $1\le i<j\le 4$, we have
\[
\rho^{(2,1^2)}_{ij} = - I_3 + v^{(2,1^2)}_{ij} {v^{(2,1^2)}_{ij}}^\T,
\]
where
\begin{eqnarray*}
&&v^{(3,1)}_{12} = \left(\begin{array}{c} \sqrt{2} \\ 0\\ 0 \end{array}\right),
v^{(3,1)}_{13} = \left(\begin{array}{c} \sqrt{1/2} \\ -\sqrt{3/2}\\ 0 \end{array}\right),
v^{(3,1)}_{14} = \left(\begin{array}{c} \sqrt{1/2} \\ -\sqrt{1/6} \\ \sqrt{4/3} \end{array}\right), \\
&&v^{(3,1)}_{23} = \left(\begin{array}{c} -\sqrt{1/2} \\ -\sqrt{3/2}\\ 0 \end{array}\right),
v^{(3,1)}_{24} = \left(\begin{array}{c} -\sqrt{1/2} \\ -\sqrt{1/6} \\ \sqrt{4/3} \end{array}\right),
v^{(3,1)}_{34} = \left(\begin{array}{c} 0 \\ \sqrt{2/3}\\ \sqrt{4/3} \end{array}\right).
\end{eqnarray*}

\noindent
For $1\le i<j\le 4$, we have
$\rho^{(1^4)}_{ij}= -1$ (one-dimensional).

{\small
\bibliography{../../../bibdb}
\bibliographystyle{amsplain}
}

\end{document}